\documentclass{amsart}

\newtheorem*{theorem}{Theorem}
\newtheorem*{lemma}{Lemma}

\begin{document}

\title{The uncertainty principle for operators determined by Lie groups}
\author{Jens Gerlach Christensen}
\address{Matematisk Institut\\ K{\o}benhavns Universitet}
\curraddr{Universitetsparken 5\\
         2100 K{\o}benhavn {\O}\\
         Denmark}
\email{vepjan@math.ku.dk}
\urladdr{http://www.math.ku.dk/\textasciitilde vepjan}
\thanks{The author wishes to thank Henrik Schlichtkrull for his proofreading
        and kind supervision}
\begin{abstract}
  For unbounded operators $A,B$ and $C$ in general, the 
  commutation relation
  $[A,B]=C$
  does not lead to the uncertainty relation
  $\| Au\|\| Bu\| \geq \frac{1}{2}|\langle \overline{C} u,u\rangle |.$
  If $A,B$ and $C$ are part of the generators of
  a unitary representation of a Lie group then the 
  uncertainty principle above holds. 
\end{abstract}

\maketitle

If $A$ and $B$ are self-adjoint or skew-adjoint operators then
\begin{equation} \label{eq:1}
  \| Au \|\| Bu \| \geq \frac{1}{2}|\langle [A,B]u,u\rangle |
   \qquad\text{for all $u\in D([A,B])$}
\end{equation}

The operator $[A,B]$ need not be closed, and if 
$C=\overline{[A,B]}$, the extended inequality
\begin{equation} \label{eq:2}
  \|Ax\|\|Bx\| \geq \frac{1}{2}|\langle Cx,x\rangle |\qquad
   \text{for all $x\in D(A)\cap D(B)\cap D(C).$}
\end{equation}
need not be valid.
K. Kraus \cite{kraus} showed that 
(\ref{eq:2}) holds when 
$A,B$ and $C$ are infinitesimal generators of a unitary
representation of a Lie group of dimension $\leq 3$. 
The result of K. Kraus was slightly extended by
G.B. Folland and A. Sitaram in \cite[Theorem 2.4]{folsit},
but the necessity of the dimension constraint was left as
an open problem.  
In the present paper the dimension constraint is shown to
be unnecessary, by means of some results  
by I.E. Segal in \cite{segal}.

If $(H,\pi)$ is a unitary representation of a Lie group $G$
and $R$ is a right-invariant vector field, then
by Stone's theorem one can define the closed skew-adjoint (by Lemma 3.1.13 in 
\cite{segal}) operator $\pi(R)$ such that 
$\pi(\exp tR)=\exp(\pi(tR))$ for all $t\in\mathbb{R}.$
Then $[\pi(X),\pi(Y)]\subseteq \pi([X,Y])$.

\begin{theorem} 
  Let $G$ be a Lie group with Lie algebra $\mathfrak{g}$, and
  let $(H,\pi)$ be a unitary representation of $G$.
  Suppose that $X,Y\in\mathfrak{g}$. 
  Then 
  \begin{equation*}
    \|Ax\|\|Bx\| \geq \frac{1}{2}|\langle Cx,x\rangle |\qquad
     \text{for all $x\in D(A)\cap D(B)\cap D(C)$}
  \end{equation*}
  holds with $A=\pi(X)$,$B=\pi(Y)$ and $C=\pi([X,Y])$.
\end{theorem}

\begin{lemma}
  Let $\sigma$ be the modular function on $G$.
  Given a right-invariant vector field $R$ on $G$, let
  $L$ be the corresponding left-invariant vector field, and 
  let $\gamma = \frac{d}{dt}\sigma(\exp(tR))|_{t=0}$.
  \begin{enumerate}
  \item If $f\in C_c^\infty(G)$ and 
        $f^*(g)= \overline{f(g^{-1})}\sigma(g)$, then
        $\pi(Rf^*)^* = \pi(Lf)+\gamma\pi(f)$. \label{item1}
  \item There is a sequence $\{f_n\}$ in $C_c^\infty(G)$
        such that $\pi(f_n)u \to u$ and
        $\pi(Rf_n)u + \pi(Lf_n)u + \gamma\pi(f_n)u \to 0$ for all
        $u\in H$. \label{item2}
  \end{enumerate}
\end{lemma}

\begin{proof}
  (\ref{item1}) is simply Lemma 3.1.7 in \cite{segal} using the fact that
  $f^{**}=f$.
  (\ref{item2}) is Lemma 3.1.8 and Lemma 3.1.9 in \cite{segal}, and 
  $\pi(f_n)u\to u$ has been proved in \cite[p.56]{knapp}.
\end{proof}

To prove the theorem I will first show that for $R\in\mathfrak{g}$ 
it holds that
$\pi(R)\pi(f_n)u \to \pi(R)u$ for $u \in D(\pi(R))$.
Fix $R\in\mathfrak{g}$ and
let $L$ be the left-invariant vector field corresponding to $R$.
Note that for all $f\in C^\infty_c(G)$ it holds that
$\pi(f)\pi(R) \subseteq -\pi(Rf^*)^*,$ since
$\langle \pi(f)\pi(R)u,v\rangle = \langle u,-\pi(Rf^*)v\rangle$ 
for all $u\in D(\pi(R))$ and $v\in H.$ 
Hence 
\begin{equation*}
  \pi(f)\pi(R) \subseteq -\pi(Lf)-\gamma \pi(f)
\end{equation*}
by (\ref{item1}) of the lemma.
But then 
\begin{align*}
  \pi(R)\pi(f_n)u-\pi(f_n)\pi(R)u 
  &=\pi(Rf_n)u+\pi(Lf_n)u+\gamma\pi(f_n)u 
\end{align*}
for any $u\in D(\pi(R)).$
By (\ref{item2}) in the lemma this tends to $0$ and since
$\pi(f_n)\pi(R)u\to\pi(R)u$ it follows that
$\pi(R)\pi(f_n)u \to \pi(R)u$.

Now let $A,B$ and $C$ be as in the theorem.
For any $u\in D(A)\cap D(B)\cap D(C)$ it follows that 
$\pi(f_n)u\to u$,
$A\pi(f_n)u\to Au$,
$B\pi(f_n)u\to Bu$ and 
$C\pi(f_n)u\to Cu$.
Since the G{\aa}rding vector 
$\pi(f_n)u$ is in $D([A,B])$  by \cite[p.56]{knapp}
the inequality (\ref{eq:1}) gives (\ref{eq:2}) as desired.

\end{document}